\documentclass[12pt, a4paper, leqno]{article}
\usepackage{amssymb}
\usepackage{amsmath}
\title{
Quasi, twisted, and all that...\\ 
in Poisson geometry and Lie algebroid theory}

\author{Yvette Kosmann-Schwarzbach}

\date{}

\def\g{\mathfrak g}
\def\G{\Gamma}
\def\p{\varphi}

\def\T*{T^*M}

\def\d{\rm d}

\def\End{\rm End}

\def\w{\wedge}

\def\bra{[~,~]}
\def\C{C^{\infty}(M)}
\def\pl{\partial}

\newtheorem{theorem}{Theorem}
\newtheorem{lemma}{Lemma}
\newtheorem{definition}{Definition}
\newtheorem{proposition}{Proposition}

\numberwithin{equation}{section}
\numberwithin{theorem}{section}
\numberwithin{definition}{section}
\numberwithin{proposition}{section}
\numberwithin{corollary}{section}
\numberwithin{lemma}{section}

\begin{document}

\maketitle

\medskip

\hspace{9truecm}{\it Dedicated to Alan Weinstein} 

\medskip

\abstract{\it Motivated by questions from quantum group and field theories, 
we review structures on manifolds that are
weaker versions of Poisson structures, and variants of the notion of
Lie algebroid. We give a simple definition of the Courant algebroids
and introduce the notion of a deriving operator for the Courant bracket
of the double of a proto-bialgebroid.
We then describe and relate 
the various quasi-Poisson
structures,
which have appeared in the literature since 1991, 
and the twisted Poisson structures studied by {\v S}evera and Weinstein.}



\section*{Introduction}

In 1986, Drinfeld introduced both the {\it quasi-Hopf algebras}, that
generalize the
Hopf algebras defining quantum groups, and their semi-classical limits,
the {\it Lie
quasi-bialgebras}.
This naturally led to the notion of {\it quasi-Poisson Lie groups}
which I introduced in \cite{yksQ} \cite{yksJ}. 

A {\it quasi-Hopf algebra} is a bialgebra in which the
multiplication is associative
but the co-multiplication is only co-associative up to a defect measured by an
element $\Phi$ in the triple tensor product of the algebra.
Similarly, the definitions of the Lie
quasi-bialgebras and the quasi-Poisson Lie groups
involve a given element in $\bigwedge ^3 \g$, where $\g$ is the
underlying Lie algebra, which 
Drinfeld denoted by $\p$.
In a Lie quasi-bialgebra,
the bracket is a Lie bracket because it satisfies the Jacobi identity, 
but the compatible cobracket is not a
true Lie bracket on the dual of $\g$, 
because it only satisfies the Jacobi identity 
up to a defect
measured by the element $\p$. On a quasi-Poisson Lie group, there is a
multiplicative
bivector field, $\pi$, whose Schouten bracket, 
$[\pi,\pi]$, does not vanish, but is also expressed in terms of $\p$.
The desire to understand 
the {\it group-valued moment maps} and the {\it quasi-hamiltonian spaces} 
of Alekseev, Malkin and Meinrenken \cite{AMM} 
in terms of Poisson geometry led to
the study of the action of quasi-Poisson Lie groups 
on manifolds equipped with a bivector field \cite{AK}.
A special case of a quasi-Poisson structure 
on a Lie group occurs when the bivector vanishes
and only $\p$ remains, corresponding to a Lie quasi-bialgebra with 
a trivial cobracket. 
The {\it quasi-Poisson manifolds} studied in
\cite{AKM} are manifolds equipped with a bivector, on which such a
quasi-Poisson Lie group acts.

\medskip

Recently, 
closed $3$-form fields 
appeared in Park's work on string theory \cite{P}, and in the work 
on topological field theory of 
Klim{\v c}ik and Strobl, 
who recognized the appearance of a new geometrical structure
which they called WZW-Poisson manifolds \cite{KS}. 
They chose this name because the role of the background
$3$-form is analogous to that of the 
Wess-Zumino term  
introduced by Witten in a field theory with target a group, 
and more recently they proposed to shorten the
name to WZ-Poisson manifolds.
Shortly after these publications circulated as preprints,
{\v S}evera and Weinstein studied such structures in the framework of Courant
algebroid theory, calling them {\it Poisson
structures with a $3$-form background}. 
They are defined in terms of a bivector
field $\pi$ and a closed $3$-form, denoted by $\p$ in \cite{SW}, but
which we shall
denote by $\psi$ to avoid confusion with the above. Again $\pi$ is
not a Poisson bivector -- unless $\psi$ vanishes, in which case the
Poisson
structure with background 
reduces to a Poisson structure --, its Schouten bracket 
is the image of the $3$-form $\psi$ 
under the morphism of vector bundles defined by $\pi$, mapping forms to
vectors. 
{\v S}evera and Weinstein also called the Poisson structures with
background {\it $\psi$-Poisson 
structures}, or {\it twisted Poisson 
structures}. This last term has since been widely used 
\cite{Rlmp}\cite{severa}\cite{CX}\cite{BCWZ}, hence the word
``twisted'' in the title of this paper. It is
justified by a related
usage in the theory of ``twisted sheaves'', and we shall occasionally use
this term but we prefer Poisson
structure with background because,
in Drinfeld's theory of
Lie quasi-bialgebras, the words 
``twist''and ``twisting'' have a different and now standard meaning.
Section \ref{bivector} of this paper is a generalization  
of Drinfeld's theory
to the Lie algebroid setting.

\medskip

The theory of Lie bialgebras, on the one hand, 
is a special case of that of the Lie
bialgebroids, introduced by Mackenzie and Xu \cite{MX}.
It was shown by Roytenberg \cite{Rlmp} that the ``quasi'' variant 
of this notion is the framework in which the Poisson
structures with background appear naturally. 
Lie algebras, on the other hand, 
are a special case of the
{\it Loday algebras}.
Combining the two approaches, we encounter the {\it Courant algebroids} of
Liu, Weinstein and Xu \cite{LWX}, or rather their equivalent
definition in terms of non-skew-symmetric brackets.

\medskip

We shall present these {\it a priori} different 
notions, and shall show how they can be
related.
In Section \ref{review}, we give a brief overview of the various
theories just mentioned. In particular we define the
{\it proto-bialgebroids} and the {\it Lie quasi-bialgebroids}, which 
generalize the Lie quasi-bialgebras, as well as their duals, the {\it
  quasi-Lie bialgebroids}.
In Section \ref{CA}, we give a 
simple definition of the Courant algebroids, which we prove to be
equivalent to the usual definition \cite{LWX} \cite{Rphd} 
(Theorem \ref{thCA}).
Liu, Weinstein and Xu \cite{LWX} showed that 
the construction of the double of Lie bialgebroids can be accomplished 
in the
framework of Courant algebroid theory by introducing {\it Manin triples
  for Lie bialgebroids}. Along the lines of \cite{Rlmp}, 
we extend these considerations
to the case of proto-bialgebroids and, in particular, to both
``quasi'' cases. Thus, we study the more general {\it Manin pairs 
for Lie quasi-bialgebroids}. 
This is the subject of Section
\ref{double}, where we 
also introduce the notion of a {\it deriving operator} (in the spirit of
\cite{yksD} and \cite{yksDB}) for the double
of a proto-bialgebroid, and we prove an existence theorem
(Theorem \ref{theoremproto}). 
Section \ref{examples} is devoted to the study of examples.
The {\it twisting} 
of Lie quasi-bialgebroids by bivectors generalizes Drinfeld's
twisting of Lie quasi-bialgebras, and leads to the
consideration of the quasi-Maurer-Cartan equation, which generalizes the 
quasi-Poisson condition. One can twist 
a quasi-Lie
bialgebroid {\it with a closed $3$-form 
  background} by a bivector, 
and the Poisson condition with background 
appears as the condition
for the twisted object to remain a quasi-Lie bialgebroid.

The world of the ``quasi'' structures which we explore here is
certainly nothing but a small part of the realm of 
homotopy structures, $L_{\infty}$, $G_{\infty}$,
etc. See, in particular, \cite{V2003} and the articles of Stasheff
\cite{S}, Bangoura \cite{B3} 
and Huebschmann \cite{H2}.
We hope to show that these are interesting objects in themselves.

\subsubsection*{Acknowledgments}
It is a pleasure to thank Henrique Bursztyn, James D. Stasheff and Thomas
Strobl for their comments on an earlier version of this text.

\section{A review}\label{review}
Before we mention the global objects such as the generalizations of the
Poisson Lie groups, we shall recall their infinitesimal counterparts.

\subsection{Lie quasi-bialgebras, quasi-Lie bialgebras 
and proto-bialgebras}\label{algebras}

We shall not review all the details of the structures that are
weaker versions of the Lie bialgebra structure, 
but we need to recall
the definition of Lie quasi-bialgebras.
It is due to 
Drinfeld \cite{D}, while in \cite{yksJ}
and \cite{byks} 
the dual case, that of a quasi-Lie bialgebra, and
the more general case of proto-bialgebras (called there
``proto-Lie-bialgebras'') are treated.
A {\it proto-bialgebra} structure on a
vector space $F$ is defined by a quadruple of elements in
$\bigwedge^{\bullet}(F \oplus F^*) \simeq C^{\infty}T^*(\Pi F)$, where
$\Pi$ denotes the change of parity. We denote such a quadruple by 
$(\mu, \gamma, \p, \psi)$, with $\mu~: \bigwedge ^2 F \to F$,
$\gamma~: \bigwedge ^2 F^* \to F^*$, $\p \in \bigwedge ^3 F$, 
$\psi \in \bigwedge ^3 F^*$. This quadruple 
defines a proto-bialgebra if and only
if $\{\mu + \gamma + \p +\psi,\mu + \gamma + \p +\psi\} =0$, where
$\{~,~\}$ is the canonical Poisson bracket of the cotangent bundle
$T^*(\Pi F)$, which coincides with the {\it big bracket} of 
$\bigwedge^{\bullet}(F \oplus F^*)$ \cite{yksJ}. This condition is
equivalent to the five conditions which we shall write below in the
more general case of the proto-bialgebroids (see Section \ref{oids}).
If either $\psi$ or $\p$ vanishes, there remain only four non-trivial 
conditions.
When $\psi = 0$, the bracket is a Lie bracket, while the
cobracket only satisfies the Jacobi identity up to a term involving
$\p$, and
we call the proto-bialgebra a {\it Lie quasi-bialgebra}.
When $\p = 0$, the bracket 
only satisfies the Jacobi identity up to a term involving
$\psi$, while the cobracket is a Lie cobracket, and 
we call the proto-bialgebra a {\it quasi-Lie bialgebra}.
Clearly, the dual of a Lie quasi-bialgebra is a quasi-Lie bialgebra, and
conversely.

Drinfeld only considered the case $\psi = 0$.
In the English translation of 
\cite{D}, what
we call a Lie quasi-bialgebra in this paper was
translated as a quasi-Lie bialgebra, a term which we shall 
reserve for the object
{\it dual} to a Lie quasi-bialgebra. In fact, 
it is in the dual object, where 
$\p = 0$ and $\psi \neq 0$ that the algebra structure is only ``quasi-Lie''.
As another potential source of confusion, we mention that 
in \cite{Rphd} and \cite{Rlmp}, the element in $\bigwedge ^ 3 F^*$
that we denote by  
$\psi$ is denoted by $\p$, and vice-versa.

Any proto-bialgebra $((F,F^*), \mu, \gamma, \p, \psi)$ has a {\it
  double}
which is ${\mathfrak d} = F \oplus F^*$, with the Lie bracket, 
\begin{align*}
  [x,y] & = \mu(x,y) +i_{x \w y}\psi \ , \\
  [x,\xi] & = - [\xi,x] = - ad^{* \gamma}_{\xi} x + ad^{* \mu}_{x} \xi \ , \\
  [\xi, \eta] & = i_{\xi \w \eta}\p + \gamma (\xi, \eta) \ . 
\end{align*}
Here $x$ and $y \in F$, and $\xi$ and $\eta \in F^*$.

Any Lie bialgebra has, associated with it, 
a pair of Batalin-Vilkovisky algebras in duality.
The extension of this property
to Lie quasi-bialgebras, giving rise to
quasi-Batalin-Vilkovisky algebras
in the sense of 
Getzler \cite{G}, has been carried out by Bangoura \cite{B1}.
There is a notion of quasi-Gerstenhaber algebra (see \cite{Rlmp}), and
Bangoura has further 
proved that quasi-Batalin-Vilkovisky algebras
give rise to quasi-Gerstenhaber algebras \cite{B2}.
For a thorough study of these notions in the general algebraic setting,
see Huebschmann \cite{H2}.
These ``quasi'' algebras are the simplest examples of $G_{\infty}$-
and $BV_{\infty}$-algebras, in which all
the higher-order multilinear maps vanish except for the trilinear map.

\subsection{Quasi-Poisson Lie groups and moment maps with values in
homogeneous spaces}
The global object corresponding to the Lie quasi-bialgebras we have
just presented was introduced in \cite{yksJ} and called a
{\it quasi-Poisson Lie group}. It is a Lie
group with a multiplicative 
bivector, $\pi_G$, whose Schouten bracket does not vanish (so that
it is not a Poisson bivector), but is a coboundary, namely
$$
\frac{1}{2} [\pi_G,\pi_G] = \p^{L}-\p^R \ ,
$$
where $\p^L$ (resp., $\p^R$) are the left- (resp., right-)invariant
trivectors on the group with value $\p \in \bigwedge^3\g$ at the identity.
In \cite{AK}, we considered the action of a quasi-Poisson
Lie group $(G,\pi_G,\p)$ on a manifold $M$ equipped with a $G$-invariant
bivector~$\pi$. When the Schouten
bracket of $\pi$ satisfies the condition
\begin{equation}\label{qP}
\frac{1}{2} [\pi,\pi] = \p_M \ ,
\end{equation}
we say that $(M,\pi)$ is a {\it quasi-Poisson $G$-space}. Here $\p_M$
is the image of the element $\p$ in $\bigwedge ^3 \g$ under the
infinitesimal action of the Lie algebra $\mathfrak g$ of $G$ on $M$.
The quasi-Poisson $G$-space 
$(M,\pi)$ is called a {\it hamiltonian quasi-Poisson $G$-space} if there
exists 
a moment map
for the action of $G$ on $M$, which takes values in $D/G$, where $D$ is
the simply connected Lie group whose Lie algebra is the double $\mathfrak
d = \g \oplus \g^*$ of the Lie
quasi-bialgebra~$\g$. See \cite{AK} for the precise definitions.

Two extreme cases of this construction are of particular
interest. 
The first corresponds to the case where the Lie
quasi-bialgebra 
is actually a Lie bialgebra ($\p=0$), {\it i.e.}, the Manin pair
with a chosen isotropic complement defining the Lie quasi-bialgebra is
in fact a Manin triple.
Then $G$ is a Poisson Lie group and 
$D/G$ can be identified with
a dual group $G^*$ of $G$. The moment maps for the quasi-hamiltonian
$G$-spaces
reduce to the moment maps in the
sense of Lu \cite{Lu} that take values in the dual Poisson Lie group, $G^*$.
The second case is that of a Lie quasi-bialgebra with vanishing
cobracket ($\gamma = 0$), to be described in the next subsection.

\subsection{Quasi-Poisson manifolds and group-valued moment maps}

Assume that $G$ is a Lie group acting on a manifold $M$, and that 
$\g$ is a {\it quadratic Lie
algebra}, {\it i.e.}, a Lie algebra 
with an invariant non-degenerate symmetric bilinear form. We consider
the bilinear form
in $\g \oplus \g$ defined as the difference of the copies of 
the given bilinear form on the two terms of the direct sum.
Let $\g$ be diagonally embedded into 
$\g \oplus \g$. Then $(\g \oplus \g, \g)$ is a Manin pair, and 
we choose the anti-diagonal, $\{(x,-x)| x\in \g\}$,
as a complement of $\g \subset \g \oplus \g$. The
corresponding Lie quasi-bialgebra has vanishing cobracket, because the
bracket of two elements in the anti-diagonal is in the diagonal, and
therefore the bivector of the quasi-Poisson structure of $G$ is trivial.
With this choice of a complement, $\p$ is 
the Cartan trivector of $\g$.
In this way, we obtain the {\it quasi-Poisson
$G$-manifolds} described in \cite{AKM}.
They are pairs, $(M,\pi)$, where $\pi$ is a $G$-invariant bivector on $M$
that satisfies equation
 \eqref{qP} with $\p$ the Cartan trivector of $\g$.
The group $G$ acting on itself by means of the adjoint action is
a quasi-Poisson
$G$-manifold, and so are its conjugacy classes. The bivector $\pi_G$
on $G$ is $\sum_a e ^R_a \w e ^L_a$, where $e_a$ is an orthonormal
basis of $\g$.
Because the homogeneous space $D/G$ of
the general theory is the group $G$ itself in this case, 
the moment maps for the 
{\it hamiltonian quasi-Poisson manifolds} are group-valued.
Those hamiltonian quasi-Poisson manifolds for which
the bivector $\pi$ satisfies a non-degeneracy condition
are precisely the {\it quasi-hamiltonian manifolds} of Alekseev,
Malkin and Meinrenken~\cite{AMM}.

\subsection{Lie bialgebroids and their doubles}
Lie bialgebroids were first defined by Mackenzie and Xu \cite{MX}. 
We state the 
definition as we reformulated it in \cite{yksG}.
To each {\it Lie algebroid} $A$ are associated
\begin{itemize}
\item{a Gerstenhaber bracket, $\bra_A$, on $\G (\bigwedge ^{\bullet} A)$,}
\item{a differential, $d_A$, on $\G (\bigwedge ^{\bullet} A^*)$.}
\end{itemize}
A {\it Lie bialgebroid} is a pair, $(A,A^*)$, of Lie algebroids in duality
such that 
$d_{A^*}$ is a derivation of $\bra_A$, or, equivalently, 
$d_A$ is a derivation of $\bra_{A^*}$.

Extending 
the construction of the Drinfeld double of a Lie
bialgebra to the case of a Lie bialgebroid is a non-trivial problem,
and several solutions have been offered, by
Liu, Weinstein and Xu \cite{LWX} in terms of the Courant algebroid $A
\oplus A ^*$, 
by Mackenzie \cite{M} in terms of the double vector bundle 
$T^*A \simeq T^*A ^*$, and by Vaintrob (unpublished) and 
Roytenberg \cite{Rphd} \cite{Rlmp} in
terms of supermanifolds.
We shall describe some properties of the first and third constructions 
in Section
\ref{double}.

\subsection{Lie quasi-bialgebroids, 
quasi-Lie bialgebroids, 
proto-bialgebroids and their doubles}\label{oids}

We call attention to the fact that we shall define here both
``Lie quasi-bialgebroids'' and 
``quasi-Lie bialgebroids'' and that, as we explain below, 
these terms are not synonymous. We extend 
the notations 
of \cite{D}, \cite{yksJ}, \cite{byks} to the case of Lie algebroids.

A {\it proto-bialgebroid} $(A,A^*)$  is
defined by anchors $\rho_A$ and $\rho_{A^*}$, brackets
$\bra_A$ and $\bra_{A^*}$, and elements $\p \in 
\Gamma (\bigwedge ^3 A)$
and $\psi \in \Gamma (\bigwedge ^3 A ^*)$.
By definition, 
\begin{itemize}
\item 
The case $\psi=0$ is that of {\it Lie quasi-bialgebroids} ($A$ is a true Lie
algebroid, while $A^*$ is only ``quasi''),
\item The case $\p=0$ is that of {\it quasi-Lie bialgebroids} ($A ^*$ is a true
Lie algebroid, while $A$ is only ``quasi'').
\item The case where both
$\p$ and $\psi$ vanish is that of the {\it Lie bialgebroids}.
\end{itemize}
While the dual of a Lie bialgebroid is itself a Lie bialgebroid, the
dual of a Lie quasi-bialgebroid is a quasi-Lie bialgebroid, and
conversely.

Whenever $A$ is a vector bundle,
the space of functions on
$T^*\Pi A$, where $\Pi$ denotes the change of parity, contains the space of
sections of $\bigwedge ^{\bullet} A$, the $A$-multivectors.
In particular, the sections of $A$ can be considered as 
functions on $T^*\Pi A$. Given the canonical isomorphism,
$T^*\Pi A^* \simeq T^*\Pi A$, the same conclusion holds for the 
sections
of $\bigwedge ^{\bullet} A^*$, in particular for the sections of
$A^*$.

A {\it Lie algebroid bracket} $\bra_A$ on a vector bundle $A$
over a manifold $M$ is defined,
together with an anchor $\rho_A~: A \to TM$, by a function $\mu$ on
the supermanifold $T^* \Pi A$  (\cite{Rphd} \cite{Rlmp} \cite{Va} 
\cite{V}). Let $\{~,~\}$
denote the canonical Poisson bracket of the  
cotangent bundle. 
The bracket of two sections $x$ and $y$ 
of $A$ is the {\it derived bracket}, in the sense of \cite{yksD},
$$
[x,y]_A= \{\{x,\mu\},y\} \ ,
$$
and the anchor satisfies 
$$
\rho_A(x)f = \{ \{x, \mu\},f\} \ ,
$$
for $f \in C^{\infty}(M)$.
When $(A,A^*)$ is a pair of Lie algebroids in duality,
both $\bra_A$ together with $\rho_A$, and $\bra_{A^*}$ together with
$\rho_{A^*}$
correspond to functions, denoted
by $\mu$ and $\gamma$,
on the same supermanifold $T^*\Pi A$, taking into account the
identification of $T^*\Pi A^*$ with $T^*\Pi A$.
The three conditions in the definition 
of a {\it Lie bialgebroid} are equivalent to the 
single equation
$$
\{\mu+\gamma,\mu+\gamma\}=0 \ .
$$

More generally, the five conditions for a 
{\it proto-bialgebroid} defined by
$(\mu, \gamma, \p, \psi)$ 
are obtained from a single equation.
By definition, a {\it proto-bialgebroid}  structure on $(A,A^*)$
is a function of degree $3$ and of Poisson square $0$ on $T^*\Pi A$.
As in the
case of a proto-bialgebra, 
such a function can be written $\mu + \gamma + \p + \psi$, where $\p
\in \G \bigwedge ^3 A$  and $\psi
\in \G \bigwedge ^3 A ^*$, and satisfies
\begin{equation}\label{def}
\{\mu+\gamma+\p+\psi,\mu+\gamma+\p+\psi\}=0 \ .
\end{equation}
The definition is equivalent to the
conditions
$$ 
\left\{
\begin{array}{rcc}
\frac{1}{2} \{\mu,\mu \} + \{\gamma, \psi \} & = & 0 \ , \\

\{\mu,\gamma \} + \{\p, \psi \} & = & 0 \ , \\

\frac{1}{2} \{\gamma,\gamma \} + \{\mu, \p \} & = & 0 \ , \\

\{\mu, \psi \} & = &  0 \ , \\

\{\gamma, \p \} & = & 0 \ .
\end{array}
\right .
$$

\begin{itemize}
\item When $((A,A^*), \mu,\gamma,\p, 0)$ is a Lie quasi-bialgebroid,  
$((A,A^*), \mu,\gamma)$ is a Lie
bialgebroid if and only if 
$\{\mu, \p\}=0$.
\item 
Dually, when $((A,A^*), \mu,\gamma, 0, \psi)$ is a quasi-Lie bialgebroid,  
$((A,A^*), \mu,\gamma)$ is a Lie
bialgebroid if and only if 
$\{\gamma, \psi\}=0$.
\end{itemize}

\noindent{\bf Remark} 
In the case of a proto-bialgebra, $(F,F^*)$, the operator
$\{\mu,~.~\}$ generalizes the
Chevalley-Eilenberg coboundary operator on cochains on $F$ with values
in $\bigwedge ^{\bullet}F$. In the term $\{\mu, \p \}$, 
$\p$ should be viewed as a
$0$-cochain on $F$ with values in $\bigwedge ^3F$, and $\{\mu, \p \}$
is an element in $F^* \otimes \bigwedge ^3 F$. So is $\{ \gamma,
\gamma \}$, which is a trilinear form on $F^*$ with values in $F^*$
whose vanishing is equivalent to the Jacobi identity for $\gamma$.
In the term $\{ \mu, \psi \}$, $\psi$ should be viewed as a $3$-cochain on $F$
with scalar values, and $\{\mu, \psi \}$ is an element in 
$\bigwedge^4 F^*$. Reversing the roles of $F$ and $F^*$, one obtains
the interpretation of the other terms in the above formulas.

\subsection{Poisson structures with background (twisted Poisson structures)}
The WZW-Poisson structures introduced by Klim{\v c}ik and Strobl
\cite{KS}
were studied by {\v S}evera and Weinstein in 2001 \cite{SW}, who called
them {\it Poisson structures
  with background}, and also {\it twisted Poisson structures}.
Roytenberg has subsequently shown that they appear by a twisting 
of a quasi-Lie bialgebroid by a bivector \cite{Rlmp}. We shall review this
  approach in Section \ref{examples}. 
The integration of Poisson structures with background 
into {\it quasi-symplectic
groupoids} is the
subject of recent work of Bursztyn, Crainic, Weinstein and Zhu
  \cite{BCWZ} and of Cattaneo and Xu
  \cite{CX}.
In addition, Xu \cite{Xu} has very recently 
extended the theory of momentum maps to this setting.

\subsection{Other structures: Loday algebras, omni-Lie algebras} 
There are essentially two ways of weakening the properties of Lie algebras. 
One possibility 
is to introduce a
weakened version of the Jacobi identity, {\it e.g.}, an
identity up to homotopy: this
is the theory of $L_{\infty}$-algebras. 
The relationship of the Courant algebroids to $L_{\infty}$-algebras
was explored in \cite{RW}.

Another possibility is
to consider non-skew-symmetric brackets: this is the theory of {\it Loday
algebras}, which Loday introduced and called {\it Leibniz algebras}.
A Loday algebra is a graded vector space with a bilinear bracket of
degree $n$ satisfying the Jacobi identity,
\begin{equation}\label{jacobi}
[x,[y,z]=[[x,y],z]+(-1)^{(|x|+n)(|y|+n)}[y,[x,z]] \ ,
\end{equation}
for all elements $x$, $y$ and $z$, where $|x|$ is the degree of $x$.
In Section \ref{CA}, 
we shall describe the Loday algebra approach to Courant
algebroids, in which case there is no grading.

The ``omni-Lie algebras'' introduced by Weinstein in \cite{W}
provide
an elegant way of characterizing the Lie algebra structures on
a vector space $V$ in
terms of the graph  in $V \times {\mathfrak {gl}}(V)$ of the adjoint operator.
In the same paper, he defined the $(R,\mathcal A)$
$C$-algebras, 
the algebraic analogue of
Courant algebroids, which generalize the $(R,\mathcal
A)$ Lie algebras (also called 
Lie-Rinehart algebras or pseudo-Lie algebras), and he posed the question of
how to determine the global analogue of an ``omni-Lie algebra''.
In \cite{KW}, he and Kinyon explored this problem and  
initiated the search for the global objects 
associated to generalized Lie algebras, that would generalize Lie groups.
They proved new properties of the Loday
algebras,
showing in what sense they can be integrated to a homogeneous left
loop, {\it i.e.}, to a manifold
with a non-associative composition law, and 
they showed that the Courant brackets of the doubles
of Lie bialgebroids can be realized on the tangent spaces of 
reductive homogeneous spaces.
These global constructions are inspired by 
the correspondence between 
generalized Lie triple systems
and non-associative multiplications on 
homogeneous spaces. (Some of the results of Bertram 
\cite{Be} might prove useful in the search for global objects
integrating generalized Lie algebras.)
For recent developments, see Kinyon's lecture \cite{Ki}.

\subsection{Generalized Poisson brackets for non-holonomic mechanical systems}
Brackets of the Poisson or Dirac type that do not satisfy the Jacobi identity
appear in many geometric constructions describing
non-holonomic mechanical systems.
There is a large literature on the subject; see for instance 
\cite{ILMD}
\cite{CMR}  
and their references.
It would be very interesting to study how these constructions relate
to the various structures
which we are now considering. In his lecture \cite{M}, Marsden showed how to
state the non-holonomic equations of Lagrangian mechanics 
in terms of isotropic subbundles in the direct sum of the tangent
and cotangent bundles of the phase space, $T^*Q$, of the system under
consideration. He calls such subbundles Dirac structures on
$T^*Q$. Yet, it is only when an integrability condition is required
that these structures become examples of the Dirac structures  
to be mentioned in the next section.

\section{Courant algebroids}\label{CA}

The 
construction of the double of a Lie bialgebra with the structure of a
Lie algebra
does not extend into a construction of the double of a Lie bialgebroid
with the structure of a Lie algebroid, because the framework 
of Lie algebroid theory
is too narrow to
permit it. While it is not
the only solution available, the introduction
of the new notion of
{\it Courant algebroid} permits the solution of this problem. 

The definition of Courant algebroids, based on Courant's earlier work
 \cite{Courant},
is due to Liu, Weinstein
and Xu \cite{LWX}. It was shown by Roytenberg \cite{Rphd}
that a Courant algebroid can be equivalently defined as a 
vector bundle $E \to M$ 
with a Loday bracket
on $\Gamma E$, an anchor
$\rho~: E \to TM$ and 
a field of non-degenerate, symmetric bilinear forms $(~|~)$ on the
 fibers of $E$, 
related by a set 
of four additional properties. It was further observed by Uchino \cite{uchino}
and by Grabowski and Marmo \cite{GM} 
that the number of independent
conditions can be reduced. We now show that it can be reduced
to two properties which are 
very natural generalizations of those of a quadratic Lie
algebra.
In fact, (i) and (ii) below are generalizations to algebroids of
the skew-symmetry of the Lie bracket, and of the condition of ${\rm
  ad}$-invariance for a bilinear form on a Lie algebra, 
respectively.

\begin{definition}
A {\em Courant algebroid} is a 
vector bundle $E \to M$ 
with 
a Loday bracket on $\Gamma E$, 
{\it i.e.}, an $\mathbb R$-bilinear map satisfying the 
Jacobi identity,
$$
[x,[y,z]]=[[x,y],z] + [y,[x,z]] \ ,
$$
for all $x, y, z \in \Gamma E$, 
an anchor, $\rho: E \to TM$, which is a morphism of vector bundles,
and  
a field of non-degenerate symmetric bilinear forms $(~|~)$ on the
fibers of $E$, 
satisfying 
\begin{align*}
{\rm (i)} \quad \quad  \quad  \rho(x)(u|v) & =(x~|~[u,v]+[v,u]) \ , \\
{\rm (ii)} \quad \quad  \quad 
\rho(x)(u|v) & =([x,u]~|~v)+(u~|~[x,v]) \ ,
\end{align*}
for all $x$, $u$ and $v \in \Gamma E$.
\end{definition}

\medskip

\noindent{\bf Remark}
Property (i) is equivalent to 
$$
{\rm (i')} \quad \quad  \quad \frac{1}{2} \rho(x)(y|y)=(x|[y,y]) 
$$ 
(which is property 4 of Definition 2.6.1 in \cite{Rphd}, and property 5 of
Section 1 in
\cite{SW}).
The conjunction of properties (i) and (ii) is equivalent to property (ii)
together with  
$$
{\rm (i'')} \quad \quad  \quad ~~
(x|[y,y])= ([x,y]|y) 
$$
(which is property 5 in Appendix A in \cite{severa}).
\medskip

We now prove 
two important consequences of properties (i) and (ii) which have been
initially considered to be 
additional, independent defining properties of Courant algebroids.

\begin{theorem} \label{thCA}
In any Courant algebroid,

\noindent {\rm (iii)} the Leibniz rule is satisfied, {\it i.e.},
$$
[x,fy]= f[x,y] + (\rho(x)f) y \ ,
$$
for all $x$ and $y\in \Gamma E$
and all $f \in C^{\infty}(M)$,

\noindent {\rm (iv)} the anchor, 
$\rho$, induces a morphism of Loday algebras from $\Gamma E$ 
to $\Gamma (TM)$, {\it i.e.}, it satisfies
$$
\rho([x,y])=[\rho x,\rho y] \ ,
$$
for all $x$ and $y \in \Gamma E$.
\end{theorem}

\noindent{\it Proof} The proof of (iii), adapted from \cite{uchino}, is
obtained by evaluating $\rho(x)(fy|z)$ in two
ways. 
We first write, using the
Leibniz rule for vector fields acting on functions, 
$$
\rho(x)(fy|z)= (\rho(x)f)(y|z) +f \rho(x)(y|z) \ .
$$
Then, using property (ii) twice, we obtain
$$
([x,fy]|z) + (fy|[x,z]) = (\rho(x)f)(y|z) +f ([x,y]|z) + f (y | [x,z]) \ .
$$
and (iii) follows by the non-degeneracy of $(~|~)$.

The proof of (iv) is that of the analogous property for Lie
algebroids (see, {\it e.g.}, \cite{yksm}). It is obtained
by evaluating $[x,[y,fz]]$, for $z \in \Gamma E$, in two ways, using
both 
the Jacobi identity for the
Loday bracket $[~,~]$ and (iii). \hfill $\Box$

\medskip

It follows from the Remark together with Theorem \ref{thCA} 
and from the arguments of
Roytenberg in
\cite{Rphd} that our definition of Courant algebroids is equivalent to
that of Liu, Weinstein and Xu in \cite{LWX}.

\medskip

A {\it Dirac 
sub-bundle} (also called a {\it Dirac structure}) 
in a Courant algebroid is a maximally isotropic
sub-bundle whose space of sections is closed under the bracket.

\medskip

Courant algebroids with base a point are {\it quadratic Lie
algebras}.
More generally, Courant algebroids with a trivial anchor are
bundles of quadratic Lie algebras with a smoothly
varying structure.

The notion of a Dirac sub-bundle in a Courant algebroid with base a
point 
reduces to that of a
maximally isotropic Lie subalgebra in a quadratic Lie algebra, 
in other words, to a {\it Manin pair}.
We shall show that a Courant algebroid together with a Dirac sub-bundle
is an appropriate generalization of the notion of a Manin pair from 
the setting of Lie algebras to that of Lie algebroids.

\medskip

A deep understanding of the nature of Courant algebroids is provided
by the consideration of the {\it non-negatively 
graded manifolds}.
This notion was defined and used
by Kontsevich \cite{K}, {\v S}evera \cite{severa} (who called them 
$N$-{\it manifolds}) and T. Voronov \cite{V}. 
In \cite{R2003}, Roytenberg
showed that the non-negatively graded symplectic manifolds of degree $2$ are
the pseudo-euclidian vector bundles, and that the Courant algebroids
are defined by an additional structure, that of a 
homological vector
field, associated to a cubic hamiltonian $\Theta$ of
Poisson square $0$, preserving the symplectic structure.
The bracket and the anchor of the Courant algebroid are recovered
from this data as the derived brackets, $[x,y]= \{\{x,\Theta\},y\}$ and
$\rho(x)f= \{\{x,\Theta\},f\}$.
T. Voronov \cite{V} studied the double of the non-negatively graded $QP$-manifolds which
are a generalization of the Lie bialgebroids.

\section{The double of a proto-bialgebroid}\label{double}

We shall now explain how to generalize the construction of a double with 
a Courant algebroid structure from Lie bialgebroids to  
proto-bialgebroids.

\subsection{The double of a Lie bialgebroid}

Liu, Weinstein and Xu \cite{LWX} have shown that complementary pairs
of Dirac sub-bundles in a Courant algebroid are in one-to-one
correspondence with Lie bialgebroids:

If $E$ is a  Courant algebroid, if $E =A \oplus B$, where $A$ and $B$
are maximally isotropic sub-bundles, 
and if $\Gamma A$ and $\Gamma B$ are closed under the
bracket, then 
\begin{itemize}
\item{$A$ and $B$ are in duality, $B \simeq A ^*$,}
\item{the bracket of $E$ induces Lie algebroid brackets on $A$ and 
$B \simeq A ^*$,
with respective anchors the restrictions of the anchor of $E$ to $A$
and $A ^*$,}
\item{the pair $(A,A ^*)$ is a Lie bialgebroid.}
\end{itemize}

Conversely, if $(A,A ^*)$ is a Lie bialgebroid, the direct sum 
$A \oplus A ^*$ 
is equipped with a Courant algebroid structure such that $A$
 and $A^*$ are maximally isotropic sub-bundles, and  
$\Gamma A$ and $\Gamma (A ^*)$ are closed under the
bracket,
the bilinear form being the canonical one, defined by
$$
(x+\xi|y+\eta) = \, <\xi,y>+<\eta,x> \ ,
$$
for $x$ and $y \in \G A$, $\xi$ and $\eta \in \G (A^*)$.

\subsection{The case of proto-bialgebroids}\label{PB}

The construction which we just recalled can be extended to the
proto-bialgebroids \cite{Rlmp}. 
Let  $A$ be a vector bundle.
Recall that a proto-bialgebroid structure on $(A,A^*)$
is a function of degree $3$ and of Poisson square $0$ on $T^*\Pi A$,
that can be written $\mu + \gamma + \p + \psi$, where $\p
\in \G (\bigwedge ^3 A)$  and $\psi
\in \G (\bigwedge ^3 A ^*)$, and $\mu$ (resp., $\gamma$) defines a
bracket and anchor on $A$ (resp., $A ^*$).

The Courant bracket 
of the double, $A\oplus A ^*$, of a  
proto-bialgebroid, $(A,A ^*)$, defined by $(\mu,\gamma, \p,,\psi)$,
is the {\it derived bracket},
$$
[x+\xi,y+\eta]
=\{\{x+\xi,
\mu + \gamma + \p + \psi\},  y +\eta\} \ .
$$
Here $x$ and $y$ are sections of $A$, $\xi$ and $\eta$ are
sections of $A^*$, and $[x+\xi,y+\eta]$ is a section of $A\oplus
A ^*$. 
(The right-hand side makes sense more generally 
when $x$ and $y$ are $A$-multivectors, and
$\xi$ and $\eta$ are $A^*$-multivectors, but the resulting quantity
is not
necessarily
a section of $\bigwedge ^{\bullet}(A \oplus A^*)$.)

The anchor is defined by
$$
\{\{x+\xi,\mu+\gamma+\p+\psi \},f \} = \{\{x,\mu \},f \}
+\{\{\xi,\gamma \},f \} = (\rho_A(x)+\rho_{A ^*}(\xi))(f) \ ,
$$
for $f \in \C$.
We set $[x,y]_{\mu}= \{\{x,\mu \},y \}$ and 
$[\xi,\eta]_{\gamma}= \{\{\xi,\gamma \},\eta \}$.
The associated {\it quasi-differentials}, $d_{\mu}$ and $d_{\gamma}$, 
on $\Gamma (\bigwedge ^{\bullet} A ^*)$ and
$\Gamma (\bigwedge ^{\bullet} A )$ are
$$d_{\mu} = \{\mu, ~ \cdot ~ \} \quad \quad {\rm and}\quad \quad  
d_{\gamma}=\{ \gamma, ~ \cdot ~ \} \ ,
$$
which satisfy
$$ 
(d_{\mu})^2 + \{d_{\gamma} \psi, ~ \cdot ~ \} = 0 \ ,
\quad \quad 
(d_{\gamma})^2 + \{d_{\mu} \p, ~ \cdot ~ \} = 0 \ .
$$

We denote the interior
product of a form $\alpha$ by a multivector $x$ by $i_x\alpha$, with
the sign convention,
$$
i_{x \wedge y} =i_x \circ i_y \ ,
$$
and we
use an analogous notation for the interior product of a multivector by a form.
The Lie derivations are defined by
$L^{\mu}_x=[i_x,d_\mu]$ and $L^{\gamma}_{\xi}=[i_{\xi},d_{\gamma}]$.
We find,
for $x$ and $y \in \Gamma A$, $\xi$ and $\eta \in \Gamma (A^*)$, 
\begin{align}
[x,y]   & =  [x,y]_{\mu} + i_{x \wedge y} \psi \ , \label{1} \\
[x,\xi]  & = -i_{\xi}d_{\gamma}x + L^{\mu}_x \xi \ , \label{2} \\ 
[\xi,x] & = L^{\gamma}_{\xi} x - i_x d_{\mu}\xi \ , \label{3} \\
[\xi,\eta] & =  i_{\xi \wedge \eta} \p + [\xi, \eta]_{\gamma} \ , \label{4}
\end{align}
that is,
$$
[x+\xi, y+\eta] = [x,y]_{\mu}
 + L^{\gamma}_{\xi}y - i_{\eta}d_{\gamma}x  
+ i_{\xi \w \eta} \p + [\xi, \eta]_{\gamma} + L^{\mu}_x \eta -
 i_yd_{\mu}\xi + i_{x \w y}\psi \ .
$$
These formulas extend both the Lie bracket of the Drinfeld double of
a proto-bialgebra \cite{byks}, recalled in Section \ref{algebras},
and the Courant bracket of the double
of a Lie bialgebroid \cite{LWX}.

\subsection{Deriving operators}
If $\xi$ is a section of $A^*$, by $e_{\xi}$ we denote the operation of 
exterior multiplication by $\xi$ on $\Gamma(\bigwedge^{\bullet}A^*)$. 
In this subsection, the square brackets $\bra$ without a
subscript denote the graded commutators of endomorphisms of 
$\Gamma (\bigwedge^{\bullet} A ^*)$.

\begin{definition} We say that a differential operator 
$\mathcal D$ on $\Gamma (\bigwedge ^{\bullet} A ^*)$ is 
a {\em deriving operator} for the Courant bracket of $A \oplus A ^*$
if it satisfies
the following relations,
\begin{align}
[[i_x, \mathcal D], i_y]
& = i_{[x,y]_{\mu}} +e_{i_{x \wedge y} \psi} \ , \label{a}
\\
[[i_x, \mathcal D],  e_{\xi}] & = - i_{i_{\xi}d_{\gamma}x} +
e_{L^{\mu}_x \xi} \ , \label{b} \\
[[e_{\xi}, \mathcal D], i_x ] & = i_{L^{\gamma}_{\xi} x} - e_{i_x
  d_{\mu}\xi} \ , \label{c}
\\
[[e_{\xi}, \mathcal D], e_{\eta}] & = i_{i_{\xi \wedge \eta} \p} 
 + e_{[\xi, \eta]_{\gamma}} \ . \label{d}
\end{align}
\end{definition}

If we identify $x \in \Gamma A$ with $i_x \in \End(\Gamma (\bigwedge ^
{\bullet}A ^*))$, 
and $\xi \in \Gamma (A ^*)$ with 
$e_{\xi} \in \End(\Gamma (\bigwedge ^ {\bullet}A ^*))$, the preceding
relations become 
\begin{align*}
[[x, \mathcal D], y]
& = {[x,y]_{\mu}} +{i_{x \wedge y} \psi} \ , \\
[[x, \mathcal D], {\xi}] & = - {i_{\xi}d_{\gamma}x} + {L^{\mu}_x \xi} \ ,
\\
[[{\xi}, \mathcal D], x ] & = {L^{{\gamma}}_{\xi} x} - {i_x d_{\mu}\xi} \ ,
\\
[[{\xi}, \mathcal D], {\eta}] & = {i_{\xi \wedge \eta} \p} 
 + {[\xi, \eta]_{\gamma}} \ ,
\end{align*}
so that the Courant bracket defined in Section \ref{PB} can also be 
written as a
derived bracket \cite{yksD} \cite{yksDB}.

\medskip

\noindent{\bf Remark}
With the preceding identification, 
the relation $i_x e_{\xi}  + e_{\xi} i_x = \, <\xi,x>$ implies that 
$$
(x+\xi)(y+\eta)+(y+\eta)(x+\xi)= (x+\xi|y+\eta) \ .
$$
This shows that $\G (\bigwedge ^{\bullet} A^*)$ is a Clifford module
of the Clifford bundle of $A \oplus A^*$, the point of departure of 
Alekseev and Xu in \cite{AX}.
 
\medskip

Does the Courant bracket of a proto-bialgebroid admit a deriving operator?
We first treat the case of a Lie bialgebroid.
The space $\Gamma (\bigwedge
^{\bullet}A ^*)$ has the structure of a Gerstenhaber algebra 
defined by $\gamma$. We shall assume that this Gerstenhaber algebra
admits a {\it generator} in the following sense \cite{Ko}. 

\begin{definition}
Let $\bra_{\mathcal A}$ be any Gerstenhaber bracket on an associative,
graded commutative algebra
$(\mathcal A, \w)$. An operator, $\pl$, on $\mathcal A$ is a 
generator of the bracket if 
$$
[u,v]_{\mathcal A}= (-1)^{|u|}(\pl(u \wedge v) - \pl u \wedge v - 
(-1)^{|u|} u \wedge \pl v) \ ,
$$
for all $u$ and $v \in {\mathcal A}$. In particular, a {\it Batalin-Vilkovisky
algebra} is a Gerstenhaber algebra which admits a generator of square $0$.
\end{definition}

\begin{lemma}
If 
$\pl$ is a generator of bracket $[~,~]_{\mathcal A}$, then,
for all $u$ and $v \in \mathcal A$,
\begin{equation}\label{pregenerator}
[e_u,\pl]= e_{\pl u} - [u,~.~]_{\mathcal A} \ ,
\end{equation}
and 
\begin{equation}\label{generator}
[[e_u,\pl],e_v]= - e_{[u,v]_{\mathcal A}} \ ,
\end{equation}
where $e_u$ is left
$\w$-multiplication by $u \in \mathcal A$.
\end{lemma}

\noindent{\it Proof}
The first relation follows from the definitions by a short computation,
and the second is a consequence of the first, since
$$
[[e_{u},\pl],e_{v}]=
[e_{\pl u} -[u, ~ \cdot ~ ]_{\mathcal A},e_{v}] 
= -  [[u, ~ \cdot ~ ]_{\mathcal A},e_{v}] 
= - e_{[u,v]_{\mathcal A}} \ ,
$$
for all $u$ and $v \in {\mathcal A}$. \hfill $\Box$
\begin{theorem}\label{theorem}
If $\pl_*$ is 
a generator of the Gerstenhaber bracket of $\Gamma (\bigwedge
^{\bullet}A ^*)$, then $d_{\mu} - \pl_*$ is a
deriving operator for the Courant bracket of 
$A \oplus A ^*$.    
\end{theorem}
 
\noindent{\it Proof}
We consider various operators acting on sections
of $\bigwedge ^{\bullet}A ^*$.
We recall from \cite{banach} (see \cite{Ko} for the case $A = TM$) 
that, for any $x
\in \Gamma A$,
\begin{equation}
[i_x,\pl_*]=-i_{d_{\gamma}x} \ .
\end{equation}
We shall also make use of the following relations,
\begin{equation}\label{comm}
[e_{\xi},d_{\mu}] = e_{d_{\mu}\xi} \ ,
\end{equation}
for $\xi \in A ^*$, 
\begin{equation}
[i_u,e_{\xi}]=(-1)^{|u|+1} 
i_{i_{\xi}u} \ ,
\end{equation}
for any $u \in \Gamma (\bigwedge ^{\bullet}A)$, and
\begin{equation}
[e_{\zeta},i_x]=(-1)^{|\zeta|+1} e_{i_x\zeta} \ ,
\end{equation}
for all $\zeta \in
\Gamma (\bigwedge ^{\bullet}A ^*)$.

\medskip

1) Let $x$ and $y$ be in $\Gamma A$. We compute 
$[i_x, d_{\mu} -\pl_*] = L^{\mu}_x + i_{d_{\gamma}x}$, whence
$$[[i_x, d_{\mu} -\pl_*],i_y] = i_{[x,y]_{\mu}} +[i_{d_{\gamma}x}, i_y]= i_{[x,y]_{\mu}} \ . $$
This proves \eqref{a}, corresponding to \eqref{1}.

2) Let $x$ be in $\Gamma A$ and let $\xi$ be in $\Gamma (A ^*)$.
We compute 
$$[[i_x, d_{\mu} -\pl_*],e_{\xi}] = [L^{\mu}_x, e_{\xi}] + [i_{d_{\gamma}x}, e_{\xi}]
= e_{L^{\mu}_x\xi} -i_{i_{\xi}d_{\gamma}x} \ . $$
This proves \eqref{b}, corresponding to \eqref{2}.

3) Since $\pl_*$ is a generating operator of $\bra_{\gamma}$,  
\eqref{pregenerator} is valid and therefore
\begin{equation}\label{cdot}
[e_{\xi}, \pl_*] = e_{\pl_*\xi} -[\xi,~ \cdot ~]_{\gamma} \ .
\end{equation}
Since $\pl_* \xi$ is of degree $0$, $e_{\pl_* \xi}$ commutes with
  $i_x$.
Therefore
$$
[[e_{\xi}, d_{\mu}-\pl_*],i_x]= [e_{d_{\mu}\xi},i_x] -
[e_{\pl_*\xi},i_x] +[[\xi,~ \cdot~ ]_{\gamma},i_x]  =
-e_{i_xd_{\mu}\xi}+ [[\xi,~ \cdot ~]_{\gamma},i_x] \ .
$$
Let us now prove that the derivation $[[\xi, ~\cdot~ ]_{\gamma},i_x]$ of 
$\Gamma (\bigwedge ^{\bullet}A ^*)$ coincides with the derivation
$i_{L^{\gamma}_{\xi}x}$.
In fact, they both vanish on $0$-forms, and on a
$1$-form $\alpha$,
$$
[\xi,i_x\alpha]_{\gamma} -i_x[\xi,\alpha]_{\gamma}= \:
\rho_{A^*}(\xi)<\alpha,x> - <[\xi,\alpha]_{\gamma}, x> \ ,
$$
while
$$
i_{L^{\gamma}_{\xi}x} (\alpha)=
\: <\alpha,L^{\gamma}_{\xi}x> \: = \: <\alpha,i_{\xi}d_{\gamma}x> +<\alpha, d_{\gamma}<\xi,x>>
$$
$$
=d_{\gamma}x(\xi,\alpha) +d_{\gamma}<\xi,x>(\alpha)
=\rho_{A^*}(\xi)<\alpha,x> - <[\xi, \alpha]_{\gamma}, x> \ .
$$
Thus
$$
[[e_{\xi}, d_{\mu}-\pl_*],i_x]=-e_{i_xd_{\mu}\xi}+i_{L^{\gamma}_{\xi}x} \ ,
$$
and \eqref{c}, corresponding to \eqref{3}, is proved.

4) Let $\xi$ and $\eta$ be sections of $A ^*$.
Then
$$
[[e_{\xi},d_{\mu}],e_{\eta}]=[e_{d_{\mu}\xi},e_{\eta}]=0 \ ,
$$
while, by \eqref{generator},
$$
[[e_{\xi},\pl_*],e_{\eta}] = - e_{[\xi,\eta]_{\gamma}} \ ,
$$
proving \eqref{d}, corresponding to \eqref{4}. \hfill $\Box$

\medskip

We now turn to the case of a proto-bialgebroid, defined by 
$\p \in \Gamma (\bigwedge ^3 A)$ and 
$\psi \in \Gamma (\bigwedge ^3 A ^*)$. 
The additional terms in the four expressions to be evaluated are

\medskip

1) $[[i_x,i_{\p}],i_y] = 0$, and 
$[[i_x,e_{\psi}],i_y] = [e_{i_x\psi},i_y]= e_{i_{x \w y}\psi}$.

\medskip

2) $[[i_x,i_{\p}],e_{\xi}] = 0$, and 
$[[i_x,e_{\psi}],e_{\xi}] = [e_{i_x\psi},e_{\xi}]= 0$.

\medskip

3) $[[e_{\xi}, i_{\p}],i_x] = [i_{i_{\xi}\p},i_x] = 0$, and 
$[[e_{\xi}, e_{\psi}],i_x] = 0$.

\medskip

4) $[[e_{\xi}, i_{\p}],e_{\eta}] = [i_{i_{\xi}\p}, e_{\eta}]
= i_{i_{\xi \w \eta} \p}$, and 
$[[e_{\xi}, e_{\psi}],e_{\eta}] = 0$.

\medskip

\noindent{Therefore, we can generalize Theorem \ref{theorem}  
as follows.}

\begin{theorem}\label{theoremproto}
If $\pl_*$ is a generator of the Gerstenhaber bracket of $\Gamma (\bigwedge
^{\bullet}A ^*)$, then $d_{\mu} - \pl_*+ i_{\p}+e_{\psi}$ is a
deriving operator for the Courant bracket of the double, 
$A \oplus A ^*$, of the proto-bialgebroid $(A,A ^*)$
defined by $\p \in \Gamma (\bigwedge ^3 A)$ and 
$\psi \in \Gamma (\bigwedge ^3 A ^*)$. 
\end{theorem}

It is clear that the addition to a deriving operator 
of derivations $i_{x_0}$ and $e_{\xi_0}$ 
of the associative, graded commutative  algebra
$\Gamma(\bigwedge^{\bullet}A ^*)$ 
will furnish a new deriving operator.
The importance of the notion of a deriving operator comes from the
fact that, if we can modifiy $d_{\mu}$ and $\partial_*$
by derivations of $\Gamma (\bigwedge
^{\bullet}A ^*)$ in such a way that the deriving operator has square
$0$, 
then the Jacobi identity for the resulting non-skew-symmetric
bracket follows
from the general properties of derived brackets that were proved in
\cite{yksD}.

Let $(A, \mu)$ be a Lie algebroid, let 
$(A,A^*)$ be the triangular Lie bialgebroid defined by a 
bivector $\pi \in \Gamma ( \bigwedge ^2 A)$ satisfying
$[\pi,\pi]_{\mu}=0$, and let 
$d_{\pi} = [\pi, ~.~]_{\mu}$ be the differential on
$\Gamma(\bigwedge^{\bullet}A)$ (see Section \ref{triang} below).
We assume that there exists a nowhere 
vanishing section, $\nu$, of the top exterior power of the dual.
Let $\pl_{\nu}$ be the generator 
of the Gerstenhaber bracket of $\Gamma(\bigwedge^{\bullet}A)$ defined by
$\nu$, which is a generator of square $0$. We set 
$$
x_{\nu}= \pl_{\nu}\pi \ .
$$
Then, $x_{\nu}$ is a section of $A$, which is called the {\it
modular field} of $(A,A^*)$ associated with~$\nu$
\cite{banach}. 
We shall now give a short proof of the existence of a deriving operator of
square~$0$ for the Courant bracket of the dual of
$(A,A ^*)$.

\begin{theorem}\label{dertriang}
The operator $d_{\pi} - \pl_{\nu} + e_{x_{\nu}}$ is a deriving operator of
square $0$ of the Courant bracket of the double of the 
Lie bialgebroid $(A^*,A)$.
\end{theorem}

\noindent {\it Proof}
By definition, the Laplacian of the strong differential
Batalin-Vilkovisky algebra 
$(\Gamma(\bigwedge^{\bullet}A), \pl_{\nu}, d_{\pi})$ is 
$[d_{\pi}, \pl_{\nu}]$,  
and we know that it satisfies the relation
$$
[d_{\pi}, \pl_{\nu}] = L^{\mu}_{x_{\nu}} \ .
$$
(See \cite{banach}, and \cite{Ko} for the case of a Poisson manifold.) 
Since, by Theorem \ref{theoremproto}, the operator $d_{\pi} - \pl_{\nu}$ is a
deriving operator, and since this property is not modified by the
addition of the derivation $e_{x_{\nu}}$, it is enough to prove that the
operator $d_{\pi} - \pl_{\nu} + e_{x_{\nu}}$ 
is of square~$0$. In fact, since both $d_{\pi}$ and
$\pl_{\nu}$ are of square $0$,
$$
\frac{1}{2} [d_{\pi}-\pl_{\nu},d_{\pi}-\pl_{\nu}] = - L^{\mu}_{x_{\nu}} =
- [x_{\nu}, ~\cdot ~]_{\mu} \ .
$$
Therefore 
$$
\frac{1}{2}
[d_{\pi}-\pl_{\nu}+e_{x_{\nu}},d_{\pi}-\pl_{\nu}+e_{x_{\nu}}] = 
- [x_{\nu}, ~\cdot~]_{\mu} + [e_{x_{\nu}},d_{\pi}] -
[e_{x_{\nu}},\pl_{\nu}] \ .
$$
By \eqref{comm}, $[e_{x_{\nu}},d_{\pi}]= e_{d_{\pi}x_{\nu}}$, which
  vanishes since $x_{\nu}$ leaves $\pi$ invariant, while
by \eqref{cdot}, $[e_{x_{\nu}},\pl_{\nu}] = e_{\pl_{\nu}x_{\nu}}-
[x_{\nu}, ~\cdot ~]_{\mu}$. In addition,
$\pl_{\nu}x_{\nu} = 0$, since $x_{\nu}=\pl_{\nu}\pi$ and
$\pl_{\nu}$ is of square $0$.
Therefore the square of $d_{\pi} - \pl_{\nu} +
e_{x_{\nu}}$ vanishes. \hfill $\Box$

\medskip

In particular, if $(M,\pi)$ is a Poisson manifold,
we obtain a deriving operator of square $0$ 
of the Courant algebroid, double of the Lie bialgebroid $(T^*M,TM)$,
dual to the triangular Lie bialgebroid $(TM,T^*M)$.

\bigskip

More generally, Alekseev and Xu \cite{AX} consider deriving operators
of the Courant bracket of a Courant algebroid whose
square is a scalar function, which they call ``generating operators''
(but which should not be confused 
with the generating operators of Batalin-Vilkovisky algebras).
They show that there always exists such a generating operator 
for the double of a Lie bialgebroid, $(A,A^*)$,
and that its square is expressible in terms of 
the modular fields of $A$ and $A ^*$ (see Theorem 5.1 and 
Corollary 5.9 of \cite{AX}). It is 
easily seen that the case of a 
triangular Lie bialgebroid is 
a particular case of their theorem and corollary,
in which the generating operator is equal to 
the deriving operator of Theorem \ref{dertriang}, and the square 
of the generating operator actually vanishes. In fact, in
the case of a 
triangular Lie bialgebroid $(A,A^*)$, 
the Laplacian $[d_{\mu},\pl_{\pi}]$
of the
strong differential Batalin-Vilkovisky algebra 
$(\Gamma(\bigwedge^{\bullet}A^*),
\pl_{\pi}, 
d_{\mu})$ vanishes 
because $\pl_{\pi} = [i_{\pi} , d_{\mu}]$, and therefore   
the modular field of $A$ vanishes. In addition, $x_{\nu}=\frac{1}{2}X_0$,
where $X_0$ is the
modular field of $A^*$ \cite{ELW} and $\pl_{\nu}x_{\nu} = 0$.
Hence, in the expression for the square of the generating 
operator given in \cite{AX}, both terms vanish.

\section{Examples}\label{examples}
We shall first analyze various constructions of Lie bialgebroids, Lie
quasi-bialgebroids and quasi-Lie bialgebroids, then we shall
consider the
Courant brackets in the theory of Poisson structures with background.
 
\subsection{Twisting by a bivector}\label{bivector}

\subsubsection{Triangular Lie bialgebroids}\label{triang}

Let $(A,\mu)$ be a Lie algebroid, and let ${\pi}$ be a
section of $\bigwedge ^2 A$. On the one hand,
such sections generalize the $r$-{\it matrices} and {\it twists}
of Lie bialgebra theory, and on the other hand,
when $A=TM$, such sections are {\it bivector fields}
on the manifold $M$. By extension, a section of
$\bigwedge ^2 A$ is called an {\it $A$-bivector}, or simply a {\it
  bivector}. 

Let $\pi ^{\sharp}$ be the vector bundle map from
$\Gamma (A ^*)$ to $\Gamma A$ defined by $\pi ^{\sharp} (\xi) =
i_{\xi}\pi$, for $\xi \in \Gamma(A ^*)$.
Consider the bracket on $A ^*$ depending on both
$\mu$ and $\pi$ defined by
\begin{equation}\label{Koszul}
[\xi,\eta]_{\mu,\pi} 
= L^{\mu}_{\pi ^{\sharp} \xi} \eta - 
L^{\mu}_{\pi ^{\sharp} \eta} \xi -
d_{\mu}(\pi(\xi,\eta)) \ ,
\end{equation}
for $\xi$ and $\eta \in \Gamma (A^*)$.
The following relation 
generalizes the equation $\gamma = - d_{\mu} r$
which is valid in a coboundary Lie bialgebra.

\begin{theorem}
Set 
\begin{equation}\label{cob}
\gamma_{\mu, \pi} = \{ \pi, \mu \} = - \{ \mu, \pi \} \ .
\end{equation} 
Then

\noindent  
{\rm (i)} 
the associated quasi-differential on $\Gamma(\bigwedge ^{\bullet}A)$
is 
\begin{equation}\label{lichne}
d_{\pi} = [\pi, ~ \cdot ~ ]_{\mu} \ ;
\end{equation}

\noindent 
{\rm (ii)} 
bracket $[\xi,\eta]_{\mu,\pi} ,$ defined
by formula \eqref{Koszul}, is equal to the derived bracket,
$$
\{\{\xi,\gamma_{\mu,\pi}\},\eta \} \ ;
$$

\smallskip

\noindent
{\rm (iii)}
if, in addition,
\begin{equation}\label{phi}
\p = - \frac{1}{2} [{\pi},{\pi}]_{\mu} \ ,
\end{equation} 
then
$((A,A ^*), \mu, \gamma_{\mu,\pi}, \p, 0)$ is  a Lie quasi-bialgebroid.
\end{theorem}

\noindent{\it Proof} The proof of (i) is a straightforward application
of the Jacobi identity.
To prove (ii) it suffices to prove that 
the quasi-differential $d_{\pi}$ is given by the usual Cartan formula
in terms of the anchor $\pi ^{\sharp}$ and the Koszul bracket
\eqref{Koszul}.
This now classic result was first proved by 
Bhaskara
and Viswanath in \cite{BV}, in the case of a
Poisson bivector on a manifold, when $A = TM$ and $[\pi,\pi]_{\mu}=0$.
We proved it independently, and in the general case, in 
\cite{yksm}. 
To prove (iii), use the
relations 
$\{ \mu, \mu \} = 0$, and
$\{ \mu, \gamma_{\mu,\pi} \} = 0$ which follows from \eqref{cob} and
the Jacobi identity. Moreover
$\{\{\mu, \pi \}, \{\mu, \pi \} \}
= \{\mu, [\pi, \pi]_{\mu} \}$,
whence 
$\frac{1}{2} \{\gamma_{\mu,\pi}, \gamma_{\mu,\pi}
 \} + \{ \mu, \p \} = 0 $,
and 
$$
- 2 \{ \gamma_{\mu,\pi}, \p \} 
= \{ \{ \pi, \mu\},[\pi,\pi]_{\mu}\} = [\pi, [\pi,\pi]_{\mu}]_{\mu} = 0 \ .
$$ 
Thus the four conditions equivalent to \eqref{def} are
satisfied. \hfill $\Box$

\medskip

The square of $d_{\pi}$ does not vanish in general, 
$$
(d_{\pi})^2  +
[\p, ~ \cdot ~ ]_{\mu} = 0 \ .
$$

A necessary and sufficient condition for $((A,A^*), \mu,
\gamma_{\mu,\pi})$ to be a  Lie bialgebroid is the 
{\it generalized Poisson condition},
\begin{equation}\label{GYB}
d_{\mu}([\pi,\pi]_{\mu})= 0 \ ,
\end{equation}
which
includes, as a special case, the generalized classical
Yang-Baxter equation, and which is equivalent to the conditions 
to be found in
\cite{yksm}, page 74, and in Theorem 2.1 in \cite{LX}.

A sufficient condition is that 
${\pi}$ satisfy the {\it Poisson condition},
\begin{equation}\label{YB}
[{\pi},{\pi}]_{\mu} =0 \ ,
\end{equation}
which generalizes both the classical Yang-Baxter equation and the
definition of Poisson bivectors.
This condition is satisfied if and only if the graph of $\pi$ is a
Dirac sub-bundle of the standard
Courant algebroid, $A \oplus A ^*$, the double of the Lie bialgebroid
with trivial cobracket, 
$((A, A ^*), \mu, 0)$. (See \cite{Courant} for the case where
$A=TM$, and \cite{LWX}.)
The Lie bialgebroid defined by $(A,\pi)$, where $\pi$ satisfies
\eqref{YB}
is called a {\it
triangular Lie bialgebroid} \cite{LX}.

\medskip

By Theorem \ref{theoremproto}, 
a deriving operator for the Courant bracket of the double of the Lie
quasi-bialgebroid $((A,A ^*), \mu, \gamma_{\mu,\pi}, - \frac{1}{2}
[{\pi},{\pi}]_{\mu},0)$ is 
$$
d_{\mu} - \pl_{{\pi}} +i_{\p} \ ,
$$
where 
$\pl_{{\pi}}$ is the graded commutator $[i_{\pi},d_{\mu}]$, 
and $\p= -\frac{1}{2}
[\pi,\pi]_{\mu}$.
In fact \cite{Ko}
\cite{yksG}, 
$\pl_{\pi}$ generates the bracket $\bra_{\mu,\pi}$ of $A ^*$.
If $\pi$ satisfies the Poisson condition \eqref{YB}, then
$d_{\mu}-\pl_{\pi}$ is a deriving operator.
\medskip

Dually, 
$((A ^*,A), \gamma_{\mu, \pi}, {\mu}, 0, \psi)$,
with $\psi = - \frac{1}{2} [{\pi},{\pi}]_{\mu}$,
is a quasi-Lie bialgebroid,
and $((A ^*,A), \gamma_{\mu, \pi}, {\mu})$ 
is a Lie bialgebroid if and only if
${\pi}$ satisfies equation \eqref{GYB}.

\subsubsection{Twisting of a proto-bialgebroid}\label{general}

The Lie quasi-bialgebroid $(A,A^*)$ and the dual quasi-Lie bialgebroid
$(A ^*,A)$ are the result of the {\it twisting} by the bivector
$\pi$ of the Lie
bialgebroid with trivial  cobracket,
$((A,A ^*), {\mu}, 0)$. The operation of twisting, in this
general setting of the theory of Lie algebroids, was  
defined and studied by Roytenberg 
in \cite{Rlmp}. He showed that one can 
also twist a proto-bialgebroid, 
$((A,A^*), \mu,
\gamma, \p, \psi)$, 
by a bivector $\pi$.
The result is a 
proto-bialgebroid defined by $( \mu'_{\pi}, \gamma'_{\pi},\p'_{\pi},
\psi'_{\pi})$,
where 
\begin{align}
{\mu'_{\pi}} & =  {\mu} +\pi^{\sharp}\psi \ , \\
{\gamma'_{\pi}} & =  {\gamma} +\gamma_{\mu, \pi} 
+(\wedge^2 \pi^{\sharp}) \psi \ , \\
\p'_{\pi} & = \p  -d_{\gamma}\pi - \frac{1}{
2} [\pi,\pi]_{\mu} +
(\wedge^3 \pi^{\sharp}) \psi \ , \\
\psi'_{\pi} & =  \psi \ .
\end{align}
Here $\pi^{\sharp}\psi$ is the $A$-valued $2$-form on $A$ such that
$$
(\pi^{\sharp}\psi)(x,y)(\xi) = \psi (x,y, \pi^{\sharp} \xi) \ ,
$$ 
for all $\xi \in \Gamma (A^*)$,
and $(\bigwedge^2 \pi^{\sharp}) \psi$ is the $A^*$-valued $2$-form
on $A^*$ 
such
that, 
$$
((\wedge^2 \pi^{\sharp}) \psi)(\xi,\eta) (x) = \psi(\pi^{\sharp} \xi,
\pi^{\sharp} \eta, x) \ ,
$$
for all $x \in \Gamma A$,
while
$(\bigwedge^3 \pi^{\sharp}) \psi$ is the
section of $\bigwedge^3 A$ such that, for $\xi, \eta$ and $\zeta \in
\Gamma (A^*)$,
$$
((\wedge^3 \pi^{\sharp}) \psi)(\xi,\eta,\zeta)
 = \psi(\pi^{\sharp} \xi, \pi^{\sharp}\eta, \pi^{\sharp} \zeta) \ .
$$
A computation shows that the tensors 
introduced above satisfy the relations
\begin{align*}
\pi^{\sharp} \psi & = \{\pi, \psi\} \ , \\
(\wedge^2 \pi^{\sharp}) \psi & = \frac{1}{2}\{\pi, \{\pi, \psi\}\} \ , \\
(\wedge^3 \pi^{\sharp}) \psi & = \frac{1}{6}
\{\pi,  \{\pi, \{\pi, \psi\}\}\} \ .
\end{align*}
These relations are used to prove that 
$((A,A^*), \mu'_{\pi}, \gamma'_{\pi},\p'_{\pi}, \psi)$
is a proto-bialgebroid.

This proto-bialgebroid is a Lie quasi-bialgebroid
if and only if $\psi = 0$,
that is, if the initial object itself was a Lie quasi-bialgebroid. 

It is a quasi-Lie bialgebroid if and only if $\p'_{\pi}=0$, that is, 
\begin{equation}\label{GtwMC}
 \p  -d_{\gamma}\pi - \frac{1}{
2} [\pi,\pi]_{\mu} +
(\wedge^3 \pi^{\sharp}) \psi = 0 \ .
\end{equation}

\medskip

We now list the particular cases of this construction that lead to the
various integrability conditions to be found in the literature.

\medskip

\noindent{\bf (a) 
Twist of a Lie bialgebroid}: $(\mu,\gamma,0,0) \mapsto 
(\mu, \gamma + \gamma_{\mu,\pi}, - d_{\gamma}\pi - \frac{1}{2}
[{\pi},{\pi}]_{\mu},0)$.
The result is a Lie quasi-bialgebroid,
furthermore it is a Lie bialgebroid if and only if the bivector $\pi$
satisfies the {\it Maurer-Cartan equation},
\begin{equation}\label{MC}
d_{\gamma}\pi + \frac{1}{2}
[{\pi},{\pi}]_{\mu} =0 \ .
\end{equation}
This condition is satisfied if and only if the graph of $\pi$ is a
Dirac sub-bundle of the 
Courant algebroid, $A \oplus A ^*$, the double of the Lie bialgebroid 
$((A, A ^*), \mu, \gamma, 0,0)$ \cite{LWX}.
A necessary and sufficient condition for $((A,A ^*), \mu, \gamma + 
\gamma _{\mu, \pi})$ to be a Lie bialgebroid is the weaker condition,
$d_{\mu}(d_{\gamma}\pi + \frac{1}{2}
[{\pi},{\pi}]_{\mu} )=0$.

If the cobracket $\gamma$ of $(A,A ^*)$ is trivial,
to $(\mu,0,0,0)$ there corresponds the quadruple 
$(\mu, \gamma_{\mu,\pi}, - \frac{1}{2}
[{\pi},{\pi}]_{\mu},0)$: this is the case studied in Section
\ref{triang}.
We know that the result is a Lie quasi-bialgebroid, and it
is a Lie bialgebroid if and only if $\pi$ satisfies the Poisson
condition \eqref{YB}, and that $((A,A^*), \mu,
\gamma_{\mu,\pi})$ is a Lie bialgebroid if and only if the bivector $\pi$
satisfies the generalized Poisson condition \eqref{GYB}.

If the bracket $\mu$ of $(A,A ^*)$ is trivial,
to $(0,\gamma,0,0)$ there corresponds the quadruple 
$(0, \gamma, -d_{\gamma}\pi,0)$, which gives rise to a Lie bialgebroid
if and only if 
\begin{equation}
d_{\gamma}\pi=0,
\end{equation} 
which means that the bivector $\pi$ on $A$ is closed, when considered as 
a $2$-form on $A^*$.

\medskip

\noindent{\bf 
(b) Twist of a Lie quasi-bialgebroid}: $(\mu,\gamma,\p,0) \mapsto 
(\mu, \gamma + \gamma_{\mu,\pi}, \p'_{\pi}, 0)$, where $\p'_{\pi} =\p 
- d_{\gamma}\pi - \frac{1}{2}
[{\pi},{\pi}]_{\mu}$. 
The result is a Lie quasi-bialgebroid,
furthermore it is a Lie bialgebroid if and only if the bivector $\pi$
and the $3$-vector $\p$ satisfy the {\it quasi-Maurer-Cartan equation},
\begin{equation}\label{qMC}
d_{\gamma}\pi + \frac{1}{2}
[{\pi},{\pi}]_{\mu} =\p \ .
\end{equation}
A necessary and sufficient condition for 
the pair $((A,A ^*), \mu,
\gamma + \gamma _{\mu, \pi})$ to be a Lie bialgebroid is the weaker condition,
$d_{\mu}\p'_{\pi}=0$.

Assume that the 
cobracket $\gamma$ of $(A,A ^*)$ is trivial. Then, in order for $(\mu,
0,\p,0)$ to define a  Lie quasi-bialgebroid, the $3$-vector $\p$
must satisfy $\{\mu, \p\}=0$. In this case, condition \eqref{qMC}
reduces to 
\begin{equation}\label{quasiP}
\frac{1}{2}
[{\pi},{\pi}]_{\mu} =\p \ ,
\end{equation}
which is a {\it quasi-Poisson} condition, analogous to \eqref{qP}.

\medskip

\noindent{\bf (c) Twist of a quasi-Lie bialgebroid}: $(\mu,\gamma, 0, \psi) 
\mapsto 
(\mu'_{\pi}, \gamma'_{\pi} , \p'_{\pi}, \psi)$, where $\mu'_{\pi}= \mu
+ \pi ^{\sharp} \psi$, 
$\gamma'_{\pi} = \gamma + \gamma_{\mu,\pi} + (\bigwedge ^2
\pi ^{\sharp})\psi$ 
and  $\p'_{\pi} =  - d_{\gamma}\pi - \frac{1}{2}
[{\pi},{\pi}]_{\mu} +(\bigwedge ^3
\pi ^{\sharp})\psi$. 
The result is a proto-bialgebroid, furthermore 
it is a quasi-Lie bialgebroid if and only if the bivector $\pi$
and the $3$-form $\psi$ satisfy the 
{\it Maurer-Cartan
equation with background $\psi$} or {\it $\psi$-Maurer-Cartan equation},
\begin{equation}\label{twMC}
d_{\gamma}\pi + \frac{1}{2}
[{\pi},{\pi}]_{\mu} = (\wedge ^3
\pi ^{\sharp})\psi  \ .
\end{equation}

Assume that the 
cobracket $\gamma$ of $(A,A ^*)$ is trivial. Then, in order for $(\mu,
0,0,\psi)$ to define a  quasi-Lie bialgebroid, the $3$-form $\psi$
must be $d_{\mu}$-closed. In this case, condition \eqref{twMC}
reduces to the 
{\it Poisson condition with
background $\psi$} or {\it $\psi$-Poisson condition},
\begin{equation}\label{tw}
\frac{1}{2}
[{\pi},{\pi}]_{\mu} = (\wedge ^3
\pi ^{\sharp})\psi  \ ,
\end{equation}
to be found in \cite{P}, \cite{KS} and \cite{SW}.

\medskip

We shall now consider in greater detail two particular cases of the above
construction of a Lie
quasi-bialgebroid from a given 
Lie quasi-bialgebroid equipped with a bivector.

\subsubsection{Lie quasi-bialgebras and $r$-matrices} 
When the base manifold of a Lie algebroid is a point,
it reduces to a  Lie algebra, ${\mathfrak g}= (F,\mu)$. 
An element in $\bigwedge ^2F$ can be viewed as a $\bigwedge
^2F$-valued $0$-cochain on $\mathfrak g$.
The {\it triangular
$r$-matrices} are those elements $r$ in $\bigwedge^2 F$
that satisfy 
$[r,r]_{\mu}=0$. 
Let us explain why the twisting defined by a bivector generalizes the 
operation of {\it twisting} defined on Lie
bialgebras, and more generally on Lie quasi-bialgebras, by Drinfeld 
\cite{D}, and further studied in \cite{yksJ} and \cite{byks}.

In this case, 
formula \eqref{Koszul} reduces to 
\begin{equation}\label{twist}
[\xi,\eta]_{\mu,r} = - (d_{\mu}r)(\xi, \eta) \ .
\end{equation}
Here $d_{\mu}r$ is the Chevalley-Eilenberg coboundary of $r$, a 
$1$-cochain on $\mathfrak g$ with values in $\bigwedge ^2 \mathfrak
g$. 
This formula is indeed that of the cobracket on $F$, obtained by
twisting a Lie bialgebra with
vanishing cobracket by an element $r \in \bigwedge^2 F$ 
(see \cite{D} \cite{yksJ}). 
Formulas \eqref{lichne} and \eqref{phi} also reduce to
the known fomulas.

Then $((F,F^*), \mu, -d_{\mu}r)$ is a Lie bialgebra 
if and only if $d_{\mu}[r,r]_{\mu}=0$, {\it i.e.}, if and only if $r$
satisfies the {\it generalized classical Yang-Baxter equation}.
A sufficient condition is that $r$ satisfy the {\it classical Yang-Baxter
equation},
$[r,r]_{\mu}=0$, in which case $r$ is a triangular 
$r$-matrix. 

In this purely algebraic case, the Courant bracket of $F \oplus F^*$
is skew-symmetric, and therefore is a true Lie algebra bracket. It
 satisfies
$$
[x,\xi]= 
- [\xi,x] = 
-i_{\xi}d_{\gamma}x + i_x d_{\mu}\xi = 
- ad^{*{\gamma}}_{\xi} x + ad^{*\mu}_x \xi \ ,
$$
and therefore coincides with the bracket of the {\it Drinfeld double}.

A deriving operator for the Lie bracket of the Drinfeld double of a
Lie proto-bialgebra $((F,F^*), \mu, \gamma, \p, \psi)$ is 
 $d_{\mu}-\partial_{\gamma}+i_{\p}+e_{\psi}$, where $d_{\mu}$ (resp., 
$\pl_{\gamma}$) is
 the generalization of the Chevalley-Eilenberg cohomology (resp., 
homology) operator of $(F, \mu)$ (resp.,
$(F^*,\gamma)$) to the case where the bracket $\mu$ (resp., $\gamma$)
does
not necessarily satisfy the Jacobi identity.

\subsubsection{Tangent bundles and Poisson bivectors}\label{poisson}
When $A = TM$, the tangent bundle of a manifold $M$, a section
$\pi$ of
$\bigwedge ^2 A$ is a bivector field on $M$.
Let $\mu_{\rm Lie}$ be the function defining 
the Lie bracket of vector fields, and more generally the
   Schouten bracket of multivector fields. 
The associated differential is the de Rham differential
   of forms, which we denote by $\d$. In this case, we denote  
the bracket of forms, defined by formula
\eqref{Koszul} above, simply by
$\bra_{\pi}$ and the function $\gamma _{\mu, \pi}$ simply by 
$\gamma_{\pi}$. Thus 
$
((TM,T^*M), \mu_{\rm Lie}, \gamma_{\pi}, \p, 0) 
$,
with
$\p = - \frac{1}{2} [{\pi},{\pi}]$, 
is a Lie quasi-bialgebroid, and if  
$[{\pi},{\pi}] =0$, {\it i.e.}, ${\pi}$ is a Poisson bivector, then
$((TM,T^*M), \mu_{\rm Lie}, \gamma_{\pi})$ is a Lie bialgebroid. 
The bracket $\bra_{\pi}$ is then the
{\it Fuchssteiner-Magri-Morosi bracket} \cite{F} \cite{MM}, 
its extension to forms of
all degrees being the
{\it Koszul bracket} \cite{Ko}.

A deriving operator for the Courant bracket of the double, $TM \oplus
T^*M$, of the Lie bialgebroid of a Poisson manifold is 
${\d} - \pl_{{\pi}}$, where $\pl_{{\pi}}= [i_{\pi},{\d}]$ is the {\it Poisson
homology operator}, defined by Koszul and studied by Huebschmann
\cite{H}, and  often called the Koszul-Brylinski operator.
Indeed, it is
well known
that the operator
$\pl_{\pi}$ generates the 
Koszul bracket of forms. This was in fact the original definition
given by Koszul in \cite{Ko}. This deriving operator is of square $0$.

\medskip

We can also consider the dual object.
Whenever ${\pi}$ is a bivector field on $M$,
$
((T^*M,TM),  \gamma_{\pi}, \mu_{\rm Lie}, 0, \psi ) 
$,
with
$\psi = - \frac{1}{2} [{\pi},{\pi}]$,
is a quasi-Lie bialgebroid, which,
when ${\pi}$ is a Poisson bivector, is the Lie bialgebroid dual to
$(TM,T^*M)$.

If $M$ is orientable with volume form $\nu$,
a deriving operator for the Courant bracket of the double, $T^*M \oplus
TM$, is 
$d_{\pi} - \pl_{\nu}$, where $\pl_{\nu}= - *^{-1}d~*$ (here, $*$ is the
operator on forms defined by $\nu$). In fact, the operator
$\pl_{\nu}$ generates the 
Schouten bracket of multivector fields \cite{Ko} \cite{banach}.
To obtain a deriving operator of square $0$, we must add to $\d -
\pl_{\nu}$ the derivation $e_{X_{\nu}}$, where $X_{\nu}$ is the
modular vector field of the Poisson manifold $(M, \pi)$ associated with
the volume form $\nu$.
In the non-orientable case, one should introduce densities as in
\cite{ELW}.
If $\pi$ is invertible, with inverse $\Omega$,
then $\pl = [i_{\Omega}, d_{\pi}]$ generates the
Schouten bracket \cite{yksG} and therefore $d_{\pi} - \pl$ is a
deriving operator of square $0$ for the Courant bracket of $T^*M \oplus TM$.
 
\subsection{The Courant bracket of Poisson structures with background}
\label{last}

\subsubsection{The Courant bracket with background}\label{CBB} 
Let $(A, \mu)$ be a Lie algebroid and let $\psi$ be a 
 $3$-form on
$A$, a section of $\bigwedge ^3A ^*$. 
Then, as we remarked in Section \ref{general},
$((A, A^*), \mu, 0, 0, \psi)$,
is a quasi-Lie bialgebroid if and only if the $3$-form $\psi$ is
$d_{\mu}$-closed,
$$
d_{\mu} \psi = 0 \ .
$$ 
This is the most general quasi-Lie bialgebroid with trivial
cobracket. By definition, the functions $\mu$ and $\psi$ satisfy  
$\{ \mu, \mu \} = 0$ and $\{\mu, \psi \}=0$, so that $\mu$
defines a Lie algebroid bracket, but we obtain 
a  Lie bialgebroid if and only if
$\psi=0$.

The bracket of the double $A \oplus A ^*$ (in the case of $TM
\oplus  T^*M$) was introduced by {\v S}evera and Weinstein \cite{SW}
who called it the modified
Courant bracket or 
the {\it Courant bracket with background} $\psi$.
This bracket satisfies
$$
[x,y]=[x,y]_{\mu} + i_{x \w y} \psi  ~~,~~
[x,\xi] = L^{\mu}_x \xi  ~~,~~
[\xi,x] = -i_xd_{\mu}\xi ~~,~~
[\xi,\eta] = 0 \ ,
$$
that is
$$
[x+\xi, y+\eta] = [x,y]_{\mu} 
+ L^{\mu}_x \eta -
 i_yd_{\mu}\xi + i_{x \w y}\psi \ .
$$

By Theorem \ref{theoremproto},   
$d_{\mu}+e_{\psi}$ is a
deriving operator of the 
Courant bracket with background $\psi$. 

In the case of a Lie algebra, $(F, \mu)$, 
$[x,\xi] = - [\xi,x]= ad^{*\mu}_x \xi$.

\medskip

\noindent {\bf Remark}
In \cite{byks}, we considered the case of the most
general Lie quasi-bialgebra with trivial cobracket.
Similarly, one can consider the Lie quasi-bialgebroids of the form
$(\mu, 0, \p,0)$, with $\{ \mu, \mu \} = 0$ and $\{\mu, \p\}=0$,
and the Courant bracket with background $\p$, a $3$-vector in this case,
$$
[x,y]=[x,y]_{\mu}  ~~,~~
[x,\xi] = L^{\mu}_x \xi  ~~,~~
[\xi,x] = -i_xd_{\mu}\xi ~~,~~
[\xi,\eta] = i_{x \w y} \p   \ ,
$$
so that
$$
[x+\xi, y+\eta] = [x,y]_{\mu}
+ i_{x \w y}\p 
+ L^{\mu}_x \eta -
 i_yd_{\mu}\xi  \ .
$$ 
This case is not dual to the preceding one.

\subsubsection{Twisting of the Courant bracket with 
  background}
Let $((A,A^*), \mu, 0,0,\psi)$ be a quasi-Lie
bialgebroid with trivial cobracket, where $\psi$ is
the background $d_{\mu}$-closed
$3$-form.
For the corresponding Courant bracket with background, 
we shall describe 
the twisting defined as above by a section $\pi$
of $\bigwedge ^2 A$. 
The twisting of this quasi-Lie bialgebroid, 
a special case of the that 
described in Section \ref{general}, 
 yields a
proto-bialgebroid whose structural elements depend on $\mu, \psi$ and
$\pi$, and which we shall denote by
$(\widetilde \mu, \widetilde \gamma, \widetilde \p,\widetilde \psi)$,
\begin{align*}
\widetilde \mu & = \mu + \pi^{\sharp}\psi \ , \\
\widetilde \gamma & =  \gamma_{\mu, \pi} +(\wedge^2\pi^{\sharp}) \psi \ , \\
\widetilde \p & = - \frac{1}{2}[\pi,\pi]_{\mu} + 
(\wedge ^3\pi^{\sharp}) \psi \ , \\
\widetilde \psi & =  \psi \ .
\end{align*}
We have seen in Section \ref{general}(c)
that the resulting twisted object is a
quasi-Lie bialgebroid if and only if
$\widetilde \p =0$, {\it i.e.}, 
$\pi$ satisfies the $\psi$-{\it Poisson condition} \eqref{tw},
\begin{equation*}
\frac{1}{2}[\pi,\pi]_{\mu} = (\wedge^3 \pi^{\sharp}) \psi  \ .
\end{equation*}
It was shown in \cite{SW} that this
condition is satisfied if and only if the graph of $\pi$ is a
Dirac sub-bundle in the Courant algebroid with background,
$A \oplus A ^*$, the double
of the quasi-Lie bialgebroid 
$((A, A ^*), \mu, 0, 0, \psi)$. 
This constitutes a generalization of the property valid in the
usual case, reviewed in
Section \ref{triang}, where $\psi =0$ and
condition \eqref{tw} 
reduces to the usual Poisson condition.

The associated derivations, on $\Gamma( \bigwedge^{\bullet} A^*)$
and on $\Gamma (\bigwedge^{\bullet} A)$, are
$$
\begin{array}{lll}
d_{\widetilde \mu} & = & d_\mu + i_{\pi^{\sharp}\psi} \ , \\
d_{\widetilde \gamma}& =
 & [\pi, ~ \cdot ~ ]_{\mu} + i_{(\bigwedge^2\pi^{\sharp})\psi} \ .
\end{array}
$$
Because $\frac{1}{2}\{\widetilde \mu,\widetilde \mu \} = -
\{\widetilde \gamma,
\widetilde \psi \}$, the derivation
$d_{\widetilde \mu}$ does not have vanishing square in general. On the
other hand,
whenever $\pi$
satisfies the $\psi$-Poisson condition,
$d_{\widetilde \gamma}$ is a true differential and $\widetilde \gamma$ defines
a true Lie bracket on $\Gamma (A^*)$,
 and a true Gerstenhaber bracket on
$\Gamma( \bigwedge^{\bullet} A^*)$, the modified Koszul bracket.
 
\medskip

We now consider the Courant bracket of the associated double,
the $\pi$-twisted Courant bracket with background
$\psi$. The mixed terms are 
$
[{x},\xi]  = 
 -i_{\xi}d_{\widetilde \gamma}{x} +L^{\widetilde \mu}_{x} \xi$ 
and $[\xi,x] =
 L^{\widetilde \gamma}_{\xi}x-i_xd_{\widetilde \mu}\xi$, therefore
\begin{align}
[{x},y]  & = 
 [{x},y]_{\mu} + (\pi^{\sharp}\psi)(x,y) \ , 
\\
[{x},\xi]  & = 
- i_{\xi}[\pi,{x}]_{\mu} - (\pi ^{\sharp}\psi)({x},\pi ^{\sharp}\xi)
+ i_x d_{\mu}\xi+ i_{x \w \pi ^{\sharp}\xi }\psi +d_{\mu}<\xi,x>  \ , 
\\
[\xi,x] & =
i_{\xi}[\pi,{x}]_{\mu} + (\pi ^{\sharp}\psi)({x},\pi ^{\sharp}\xi)
+ [\pi, <\xi,x>]_{\mu} 
- i_x d_{\mu}\xi - i_{x \w \pi ^{\sharp}\xi }\psi \ ,
\\
[\xi, \eta] & =  
[\xi,\eta]_{\mu, \pi}+ i_{\pi^{\sharp}\xi \w \pi^{\sharp}\eta }\psi \ .
\label{modified}
\end{align}

In particular, for $\psi =0$, we obtain the Courant bracket of the
double of the twist by $\pi$ of the Lie bialgebroid with trivial
cobracket,
$((A,A ^*), \mu,0)$, considered in Section \ref{triang}. 
Therefore, whenever $\pi$ satisfies the generalized Poisson
condition \eqref{GYB}, the above formulas yield the Courant bracket
of the double of the Lie bialgebroid $((A,A ^*), \mu,
\gamma_{\mu,\pi})$.
In the purely algebraic case, we recover the Drinfeld double of a
coboundary Lie bialgebra, defined by $r$, an $r$-matrix solution of the
generalized Yang-Baxter equation. 
Setting $\underline
r(\xi)=i_{\xi}r$, and using the relation 
$i_x d_{\mu}\xi = ad^{*\mu}_x \xi$, we obtain 
\begin{align*}
 [x,y]  & = [x,y]_{\mu}  \ , \\
        [x,\xi] & =  
- [\xi,x] = -\underline r (ad^{*\mu}_x \xi) + ad^{\mu}_x (\underline
r\xi) + ad^{*\mu}_x \xi \ , \\
  [\xi,\eta] & =  ad ^{*\mu}_{\underline r \xi} \eta - 
ad^{*\mu}_{\underline r \eta} \xi \ .
\end{align*} 

To conclude, we prove a property 
of the Lie bracket
defined by $\widetilde \gamma$ on $\Gamma (A ^*)$.
 
\begin{proposition}
If, $\pi$ satisfies the $\psi$-Poisson condition, 
the mapping $\pi ^{\sharp}$ is a morphism of Lie algebroids fom $A ^*$
  with the Lie bracket \eqref{modified} to $A$
with the Lie bracket $\bra_{\mu}$. 
\end{proposition}

\noindent {\it Proof} 
It is clear that the anchor of $A ^*$ is $\rho_A \circ \pi ^{\sharp}$.
To prove that $\pi ^{\sharp}$ satisfies 
\begin{equation}\label{morphism}
\pi ^{\sharp} [\xi,\eta] = [\pi ^{\sharp} \xi, \pi ^{\sharp}
\eta]_{\mu} \ ,
\end{equation}
for all $\xi$ and $\eta \in \Gamma(A ^*)$, we recall the relation, 
\begin{equation}\label{pn}
\pi ^{\sharp} [\xi,\eta]_{\mu,\pi} - [\pi ^{\sharp} \xi, \pi ^{\sharp}
\eta]_{\mu} =
\frac{1}{2}[\pi,\pi]_{\mu}(\xi,\eta) \ ,
\end{equation}
proved in \cite{yksm}.
In view of \eqref{modified}, where the bracket of $A ^*$ is expressed 
in terms
of $\bra_{\mu, \pi}$ and $\psi$, and of the equality,
$$
( (\wedge ^3 \pi ^{\sharp})\psi )(\xi, \eta) = -  \pi ^{\sharp}
(i_{\pi ^{\sharp}\xi \wedge \pi ^{\sharp} \eta}\psi) \ ,
$$
we see that, when $\pi$ satisfies the $\psi$-Poisson condition \eqref{tw}, 
equation \eqref{morphism} follows from 
equation \eqref{pn}. \hfill $\Box$

\subsubsection*{Conclusion}
In the preceding discussion, we have encountered various
weakenings and generalizations of the usual notions of Lie bialgebra,
Lie algebroid and Poisson structure  
that have appeared in the literature, starting with Drinfeld's
semi-classical limit of quasi-Hopf algebras, and up to the 
recent developments due in great part
to Alan Weinstein, his co-workers and his 
former students. We hope to have clarified the relationships and
properties of these structures.

\bigskip

\noindent{Centre de Math{\'e}matiques (U.M.R. 7640 du C.N.R.S.)}

\noindent{{\'E}cole Polytechnique}

\noindent{F-91128 Palaiseau, France}

\end{document}